\documentclass[a4paper,12pt]{amsart}
\usepackage{mathrsfs}
\usepackage{amsfonts}
\usepackage{amssymb}
\usepackage{ifthen}
\usepackage{graphicx}
\usepackage[usenames]{color}
\nonstopmode \numberwithin{equation}{section}
\newtheorem{thm}{Theorem}
\newtheorem{cor}{Corollary}
\newtheorem{lem}{Lemma}

\newtheorem{claim}{Claim}
\newtheorem{conj}{Conjecture}

\theoremstyle{definition}
\newtheorem{defn}{Definition}
\newtheorem{case}{Case}
\newtheorem{subcase}{Subcase}

\newtheorem{examp}{Example}
\newtheorem{prob}{Problem}
\newtheorem{ques}[equation]{Question}
\newtheorem{rem}{Remark}

\newcounter {own}
\def\theown {\thesection       .\arabic{own}}

\newenvironment{pf}[1][]{%
 \vskip 3mm
 \noindent
 \ifthenelse{\equal{#1}{}}%
  {{\slshape Proof. }}%
  {{\slshape #1.} }%
 }%
{\qed\bigskip}

\newcounter{alphabet}

\newcounter{minutes}
\divide\time by 60
\newcounter{hours}
\multiply\time by 60
\addtocounter{minutes}{-\time}

\newcommand{\IN}{{\mathbb N}}
\newcommand{\IC}{{\mathbb C}}
\newcommand{\ID}{{\mathbb D}}

\def\be{\begin{equation}}
\def\ee{\end{equation}}

\newcommand{\bee}{\begin{enumerate}}
\newcommand{\eee}{\end{enumerate}}

\newcommand{\blem}{\begin{lem}}
\newcommand{\elem}{\end{lem}}
\newcommand{\bthm}{\begin{thm}}
\newcommand{\ethm}{\end{thm}}
\newcommand{\bcor}{\begin{cor}}
\newcommand{\ecor}{\end{cor}}
\newcommand{\beg}{\begin{examp}}
\newcommand{\eeg}{\end{examp}}
\newcommand{\begs}{\begin{examples}}
\newcommand{\eegs}{\end{examples}}
\newcommand{\bdefe}{\begin{defn}}
\newcommand{\edefe}{\end{defn}}
\newcommand{\bprob}{\begin{prob}}
\newcommand{\eprob}{\end{prob}}
\newcommand{\bques}{\begin{ques}}
\newcommand{\eques}{\end{ques}}
\newcommand{\bei}{\begin{itemize}}
\newcommand{\eei}{\end{itemize}}
\newcommand{\bca}{\begin{case}}
\newcommand{\eca}{\end{case}}
\newcommand{\bsca}{\begin{subcase}}
\newcommand{\esca}{\end{subcase}}

\newcommand{\bcl}{\begin{claim}}
\newcommand{\ecl}{\end{claim}}

\newcommand{\bcon}{\begin{conj}}
\newcommand{\econ}{\end{conj}}
\newcommand{\bcons}{\begin{conjs}}
\newcommand{\econs}{\end{conjs}}
\newcommand{\bprop}{\begin{propo}}
\newcommand{\eprop}{\end{propo}}
\newcommand{\br}{\begin{rem}}
\newcommand{\er}{\end{rem}}
\newcommand{\brs}{\begin{rems}}
\newcommand{\ers}{\end{rems}}
\newcommand{\bo}{\begin{obser}}
\newcommand{\eo}{\end{obser}}
\newcommand{\bos}{\begin{obsers}}
\newcommand{\eos}{\end{obsers}}
\newcommand{\bpf}{\begin{pf}}
\newcommand{\epf}{\end{pf}}
\newcommand{\ba}{\begin{array}}
\newcommand{\ea}{\end{array}}
\newcommand{\beq}{\begin{eqnarray}}
\newcommand{\beqq}{\begin{eqnarray*}}
\newcommand{\eeq}{\end{eqnarray}}
\newcommand{\eeqq}{\end{eqnarray*}}

\newcommand{\ra}{\rightarrow}

\newcommand{\ds}{\displaystyle}

\begin{document}
\bibliographystyle{amsplain}
\title {On Harmonic Entire mappings}

\author[H. Deng]{Hua Deng}
\address{H. Deng, Department of Mathematics, Hebei University,
Baoding, Hebei 071002, People's Republic of China.}
\email{1120087434@qq.com}

\author[S. Ponnusamy]{Saminathan Ponnusamy}

\address{S. Ponnusamy, Department of Mathematics, Indian Institute of
Technology Madras, Chennai-600 036, India. }
\email{samy@iitm.ac.in}

\author[J. Qiao]{Jinjing Qiao${}^\dagger
$}
\address{J. Qiao, Department of Mathematics, Hebei University,
Baoding, Hebei 071002, People's Republic of China.}
\email{mathqiao@126.com}
\author[Y. Shan]{Yanan Shan}
\address{Y. Shan, Department of Mathematics, Hebei University,
Baoding, Hebei 071002, People's Republic of China.}
\email{1067382557@qq.com}

\keywords{Entire function, harmonic entire mapping, order, type, univalent, close-to-convex and convex functions\\
$
{}^\dagger$ Corresponding author}


\begin{abstract}
In this paper, we investigate properties of harmonic entire mappings.
Firstly, we give the characterizations of the order and the type for a harmonic entire mapping $f=h+\overline{g}$,
respectively, and also consider the relationship between the order and the type of $f$, $h$, and $g$. Secondly, we
investigate the harmonic mappings $f=h+\overline{g}$ such that $f^{(n_p)}=h^{(n_p)}+\overline{g^{(n_p)}}$ are
univalent in the unit disk, where $\{n_p\}_{p=1}^{\infty}$ be a strictly increasing sequence of nonnegative integers.
In terms of the sequence $\{n_p\}_{p=1}^{\infty}$, we derive
several necessary conditions for these mappings to be entire and also establish an upper bound for
the order of these mappings.
\end{abstract}


\maketitle \pagestyle{myheadings} \markboth{H. Deng, S Ponnusamy, J. Qiao and Y. Shan}{On Harmonic Entire mappings}

\section{Introduction}

There are abundant results for entire functions (i.e. functions analytic in the whole complex plane $\mathbb{C}$).
See, for example, \cite{Bo1}. The chief purpose of this paper is to establish several properties of
harmonic entire functions, i.e. complex-valued harmonic functions $f$  in $\mathbb{C}$, where $u$ and $v$ are
real-harmonic in $\IC$. 
Throughout the article we refer a harmonic function $f$ in $\IC$ as a harmonic entire mapping.
Since $\IC$ is simply connected, each such an $f$ has a canonical decomposition
$f=h+\overline{g}$, where $h$ and $g$ are analytic in $\IC$, and $g(0)=0$.

A complex-valued harmonic mapping $f$ in a simply connected domain $D$, not identically constant,
will be classified as sense-preserving in $D$ if it satisfies the
Beltrami equation of the second kind, $\overline{f}_{\overline{z}}=\omega f_{z}$, where $\omega$ is an analytic
function in $D$ with $|\omega(z)|<1$ for $z\in D$.  Since the Jacobian $J_{f}$ of $f$ is given by
$J_{f}(z)=|f_{z}(z)|^{2}-|f_{\overline{z}}(z)|^{2}$, this implies in particular that $J_{f}(z)>0$ wherever
$f_{z}(z)\neq 0$ in $D$. Several properties of planar harmonic mappings of the unit disk $\ID:=\{z:\, |z|<1\}$
may be found from \cite{Du} (see also \cite{CS} and a little survey chapter \cite{PonRasi2013}).

The paper is organized as follows.
In Section \ref{csw-sec2}, we introduce the order and the type of harmonic entire mappings,
and obtain characterizations of the order and the type for harmonic entire mappings.
In Section \ref{csw-sec3}, we discuss univalent harmonic mappings $f=h+\overline{g}$ with univalent derivatives
$f^{(n)}=h^{(n)}+\overline{g^{(n)}}$ ($n\geq 0$), where we regard the relation $f^{(0)}=h^{(0)}+\overline{g^{(0)}}
$ as $f=h+\overline{g}$. In particular, we obtain that if $f$ is harmonic in $\mathbb{D}$ and $f^{(n)}$ is univalent in $\ID$ for each $n\geq 0$,
then $f$ is a harmonic entire mapping, see Theorem \ref{thm3.1}. In fact, in order to ensure that a harmonic mapping to be entire,
it is not necessary that all $f^{(n)}$ have to be univalent in $\ID$. Let $R_n$ be the largest number with the property that
$f^{(n)}$ is univalent in  $|z|<R_n$, and  $R_f$ be the radius of convergence of the power series
representation of $f$ about the origin, i.e., $R_f=\min\{ R_h, R_g\}$, where the power series which represents $h$ and $g$ converges
for $|z|<R_h$ and $|z|<R_g$, respectively. We discuss relationships between $\{R_n\}_{n=1}^{\infty}$ and $R_f$, see
Theorems \ref{thm3.2} and \ref{thm3.3}. Finally, we consider harmonic entire mappings $f$ such that $f^{(n_p)}$ are univalent
in $\mathbb{D}$, where $\{n_p\}_{p=1}^{\infty}$  is a strictly increasing sequence of nonnegative integers, and then discuss a
necessary condition for these mappings to be entire in terms of $\{n_p\}_{p=1}^{\infty}$, see Theorem \ref{thm3.4}. We also obtain
an upper bound of the order for these mappings, see Theorems \ref{thm3.5} and \ref{thm3.6}.

\section{The order and the type of harmonic entire mappings }\label{csw-sec2}

In this section, we study the rate of growth of the maximum modulus of
a harmonic entire mapping by introducing the order and the type.


Let $ 0<R\leq \infty$ and $r<R$. If $f$ is harmonic in $|z|<R$, then we denote by $M(r,f)$
the maximum modulus of $f$ for $|z|=r<R$:
\be\label{eq2.1}
M(r,f)=\max_{0\leq\theta<2\pi}|f(re^{i\theta})|.
\ee

\bdefe\label{de2.1}
We say that a harmonic entire mapping $f$ is of order $\rho$, $0\leq\rho\leq\infty$, if
$$
\rho=\limsup_{r\rightarrow\infty}\frac{\log\log M(r,f)}{\log r}.
$$
\edefe

Obviously, every constant function has order $0$. Clearly, $f$ is of finite order $\rho$ if and only if, for every positive
$\varepsilon$ but for no negative $\varepsilon$,
\be\label{eq2.3}
M(r,f)= O(\exp(r^{\rho+\varepsilon}))\quad (r\rightarrow\infty).
\ee
That is, there is a positive number $M_1(\varepsilon)$ such that $M(r,f)\leq M_1(\varepsilon) \exp(r^{\rho+\varepsilon})$.

If the order $\rho$ of $f$ is finite and different from zero, then we define another number, called the type of $f$,
which gives a more precise description of the rate of growth of $f$. See \cite{Bo1}.

\bdefe\label{de2.2}
A harmonic entire mapping $f$ of order $\rho$ ($0<\rho<\infty$) is said to be of type $\tau$, $0\leq\tau\leq\infty$, if
\be\label{eq2.4}
\tau =\limsup_{r\rightarrow\infty} \frac{\log M(r,f)}{r^{\rho}}.
\ee
\edefe

According as $\tau=\infty$, $0<\tau<\infty$, or $\tau=0$, $f$ is said to be of maximum (or infinite), mean (or normal or finite),
or minimum  (or zero) type of order $\rho$, respectively. It is of finite type $\tau$ if and only if, for every
positive $\varepsilon$ but for no negative $\varepsilon$,
$$M(r,f)=O(\exp((\tau+\varepsilon)r^\rho)).
$$

It is often convenient to have a term to describe a harmonic entire mapping which is of order not exceeding
$\rho$ and of type not exceeding $\tau$ of order $\rho$. Although there is no standard term for this,
but we shall say in this paper that such an $f$ is of growth $(\rho, \tau)$. A harmonic entire mapping of
growth $(1, \tau)$, $\tau<\infty$, is called a harmonic  entire mapping of exponential
type or exponential type $\tau$.  Harmonic entire mappings of exponential type $\tau$ necessarily include all harmonic
entire mappings of order $1$ and type less than or equal to $\tau$, as well as all harmonic entire mappings of order less than $1$.

First, we recall the following lemma.

\blem\label{lem2.2} {\rm (}\cite[Lemma 14.1.1]{Hi}{\rm)}
The entire function
\be\label{eq2.5}
F_{\rho}(z)=\sum_{n=1}^{\infty}\left(\frac{n}{\rho e}\right)^{-\frac{n}{\rho}}z^{n}, ~{\rho}>0,
\ee
is of order $\rho$ and of type $1$ of that order. More precisely, as $r\rightarrow\infty$,
$$
M(r, F_{\rho})=F_{\rho}(r)=(2\pi)^{1/2}{\rho} r^{\rho/2}\exp(r^{\rho})\{1+o(1)\}.
$$
\elem

Now, we are ready to present the first characterization of the order of a harmonic entire mapping given by the formula
\eqref{eq2.6},  which shows that there is a close relationship between the coefficients and the order.
Our first two results are natural extension of \cite[Theorems~14.1.1 and 14.1.2]{Hi} for harmonic case. The proofs of
the next two theorems follows from the method of the proof of the analytic case, but for the sake of completeness
we include the necessary details.

\bthm\label{thm2.1}
Let $f=h+\overline{g}$ be a harmonic entire mapping, where
\be\label{DPQS-eqex1}
h(z)=\sum _{n=0}^{\infty}a_nz^n~\mbox{ and}~g(z)=\sum _{n=1}^{\infty}b_nz^n.
\ee
If the order of $f$ is $\rho$, then $\rho$ equals $\sigma$, where $\sigma$ is determined by the formula
\be\label{eq2.6}
\sigma =\limsup_{n\rightarrow\infty}\frac{n\log n}{\log[1/(|a_{n}|+|b_{n}|)]}.
\ee
\ethm
\bpf First we begin by showing that $\rho\leq\sigma$.
If $\sigma$ is finite and $\varepsilon>0$ is fixed, then there exists an $n_{\varepsilon}$ such that
$$\frac{n\log n}{\log[1/(|a_{n}|+|b_{n}|)]}<\sigma+\varepsilon,  ~
\mbox{ i.e. $|a_{n}|+|b_{n}|<n^{-\frac{n}{\sigma+\varepsilon}}$,}~\mbox{for $n>n_{\varepsilon}$.}
$$
From this we get that there exists a finite $K(\varepsilon)$ such that
$$|a_{n}|+|b_{n}|<K(\varepsilon)n^{-\frac{n}{\sigma+\varepsilon}}~\mbox{ for all $n\geq 1$.}
$$

Thus, by the definition of $F_{\rho}(z)$ we have
$$M(r,f)\leq|a_{0}|+K(\varepsilon)F_{\rho_{0}}[(\rho_{0} e)^{-\frac{1}{\rho_{0}}}r],\quad \rho_{0}=\sigma+\varepsilon.
$$
Since the order of $F_{ \rho_{0}}(z)$ equals $\rho_{0}$, we see that $ \rho\leq\sigma+\varepsilon$ for every $\varepsilon>0$.
Allowing $\varepsilon \ra 0$, we obtain that $\rho\leq\sigma$ provided  $\sigma$ is finite.
But the inequality is trivially true if $\sigma=\infty$ and thus, $\rho\leq\sigma$ holds regardless of
whether $\sigma$ is finite or $\infty$.

Next we have to prove that $\sigma\leq\rho$. The first step is to obtain explicit formulas
for the linear functionals $a_{n}$ and $b_{n}$ in terms of $f=h+\overline{g}$. For $|z|<r$, we have
$$f(z)=\frac{1}{2\pi}\int^{2\pi}_{0}P(z,re^{it})f(re^{it})\,dt,
$$
where $P(z,re^{it})$ denotes the Poisson kernel given by
$$P(z,re^{it})= \frac{r^{2}-|z|^{2}}{|re^{it}-z|^{2}}=\frac{re^{it}}{re^{it}-z}+\frac{\overline{z}}{re^{-it}-\overline{z}},
~|z|<r.
$$
Expanding into geometric series, we see that
$$P(z,re^{it}) =\sum^{\infty}_{n=0}r^{-n}e^{-int}z^{n}+\sum^{\infty}_{n=1}r^{-n}e^{int}\overline{z}^{n},
$$
so that
$$a_{n}=\frac{1}{2\pi}\int^{2\pi}_{0}r^{-n}e^{-int}f(re^{it})\,dt ~\mbox{ and }~
b_{n}=\frac{1}{2\pi}\int^{2\pi}_{0}r^{-n}e^{-int}\overline{f(re^{it})}\,dt.
$$

For any $r>0$, we easily have by \eqref{eq2.1} and the last relation that
$$|a_{n}|\leq M(r,f)r^{-n} ~\mbox{ and }~|b_{n}|\leq M(r,f)r^{-n}
$$
and therefore,
\be\label{eq2.11}
|a_{n}|+|b_{n}|\leq 2r^{-n}M(r,f).
\ee
Suppose that $\rho$ is finite. Then for any fixed $\varepsilon>0$, by \eqref{eq2.3},
there is a finite $K_{1}(\varepsilon)$ such that for all $r$,
$$M(r,f)\leq K_{1}(\varepsilon)\exp(r^{\rho+\varepsilon}),~\mbox{ i.e., }~r^{-n}M(r,f)\leq K_{1}(\varepsilon)\varphi (r),
$$
where $\varphi (r)=r^{-n} \exp(r^{\rho+\varepsilon})$ for all $r>0$. Since $r$ is at our disposal, we wish to find a suitable $r$
such that $\varphi (r)$ is minimum. To do this, for a fixed $n$, we find that
$$\min _{r}\varphi (r)= \varphi (r_0)=\left (\frac{n}{e(\rho+\varepsilon)}\right )^{-\frac{n}{\rho+\varepsilon}}
, \quad r_0=\left (\frac{n}{\rho+\varepsilon}\right )^{\frac{1}{\rho+\varepsilon}}.
$$
It follows from \eqref{eq2.11} that
$$|a_{n}|+|b_{n}|\leq 2K_{1}(\varepsilon)\left (\frac{n}{e(\rho+\varepsilon)}\right )^{-\frac{n}{\rho+\varepsilon}},
$$
and thus,
$$\sigma=\limsup_{n\rightarrow\infty}\frac{n\log n}{\log[1/(|a_{n}|+|b_{n}|)]}\leq\rho+\varepsilon
$$
for every $\varepsilon>0$. This gives $\sigma\leq\rho$ and this inequality is obviously true if $\rho=\infty$.
Hence $\sigma=\rho$ as asserted.
\epf

\br
{\em
In calculating $\sigma$, for $|a_{n}|+|b_{n}|=0$ for all large $n$, we can take the right side of \eqref{eq2.6} as $0$.
}
\er

\bthm\label{thm2.2}
Let $f=h+\overline{g}$ be a harmonic entire mapping of order $\rho$,  $0<\rho<\infty$, where $h$ and $g$ are given by
\eqref{DPQS-eqex1}. If $f$ is of finite type, then the type of $f$ is given by
\be\label{eq1.9}
\tau_0=\frac{1}{e\rho}\limsup_{n\rightarrow\infty}n(|a_{n}|+|b_{n}|)^{\frac{\rho}{n}}.
\ee
\ethm
\bpf Let the type of $f$ be $\tau$.
If $\tau_0$ given by \eqref{eq1.9} is finite, then, as in the proof of Theorem \ref{thm2.1}, given $\varepsilon >0$
there exists a $K(\varepsilon)$ such that
$$|a_{n}|+|b_{n}|\leq K(\varepsilon)\left (\frac{n}{e\rho(\tau_0+\varepsilon)}\right )^{-\frac{n}{\rho}}~\mbox{ for $n\geq 1$},
$$
which in turn implies that
$$M(r,f)\leq   |a_0| +  K(\varepsilon)F_{\rho}[(\tau_0+\varepsilon)^{\frac{1}{\rho}}r],
$$
where $F_{\rho}(z)$ is defined by \eqref{eq2.5}. Now, because $F_{\rho}(z)$ is of order $\rho$ and type $1$, we get that $\tau\leq \tau_0$.

Next we prove that $\tau_0\leq\tau$. Again, let $\tau<\infty$. Then  by using \eqref{eq2.11}, it follows that
for every $\varepsilon>0$, there is a finite $K_{1}(\varepsilon)$, such that
$$M(r,f)\leq K_{1}(\varepsilon)\exp[(\tau+\varepsilon)r^{\rho}],~\mbox{ i.e., }~r^{-n}M(r,f)\leq
K_{1}(\varepsilon)\psi (r),
$$
where $\psi (r)=r^{-n} \exp[(\tau+\varepsilon)r^{\rho}]$ for  $r>0$. Thus, by the estimate \eqref{eq2.11}, we
deduce that
$$|a_{n}|+|b_{n}|\leq 2r^{-n}M(r,f)\leq  2K_{1}(\varepsilon)\psi (r)
$$
for any $r>0$ and for any fixed $\varepsilon>0$. Now, we find an $r$ such that $\psi (r)$ is minimum.
To do this, for a fixed $n$, we find that
$$\min _{r}\psi (r)= \psi (r_1)=\left (\frac{n}{e\rho(\tau+\varepsilon)}\right )^{-\frac{n}{\rho}}
, \quad r_1=\left (\frac{n}{\rho(\tau+\varepsilon)}\right )^{\frac{1}{\rho}}.
$$
Finally, we conclude that
$$
|a_{n}|+|b_{n}|\leq 2K_{1}(\varepsilon) \left (\frac{n}{e\rho(\tau+\varepsilon)}\right )^{-\frac{n}{\rho}},
$$
and thus,
$$\tau_0=\frac{1}{e\rho}\limsup_{n\rightarrow\infty}n(|a_{n}|+|b_{n}|)^{\frac{\rho}{n}}\leq\tau+\varepsilon.
$$
Since $\varepsilon>0$ is arbitrary, we deduce that $\tau_0\leq\tau$. Hence $\tau_0=\tau$ as asserted.
\epf

As $f=h+\overline{g}$ is a harmonic entire mapping if and only if $h$ and $g$ are entire functions, it
is natural to ask the following.
\bprob
{\em
What is the relationship between the order of $f$, the order of $h$ and the order of $g$?
}
\eprob
Our next result answers this question.

\bthm\label{thm2.3}
Suppose $f=h+\overline{g}$ is a harmonic entire mapping such that  $f$, $h$ and $g$ is of
order $\rho$, $\rho_{h}$ and $\rho_{g}$, respectively. Then  $\rho=\max\{\rho_{h}, \rho_{g}\}$.

Moreover, if $0<\rho_h<\infty$ and $0<\rho_g<\infty$, let the type of $f$, $h$ and $g$  be denoted by $\tau$, $\tau_h$ and $\tau_g$, respectively. Then
the following holds:
$$
\tau= \left\{
\begin{array}{ll}
\tau_h  &\mbox{ if $\rho_g<\rho_h$,}\\
\ds \tau_g    &\mbox{ if $\rho_g>\rho_h$,}\\
\ds \max\{\tau_h, \tau_g\}  &\mbox { if $\rho_h=\rho_g$.}
\end{array} \right.
$$

\ethm\bpf
Without loss of generality, we may assume that $\rho_g\leq \rho_h<\infty $ (otherwise, we consider $\overline{f}$)
so that $\rho_h=\max\{\rho_h,\rho_g\}$.

Now, for every $\varepsilon>0$, there exists an $r_{0}>1$ such that
$$ M(r,h)\leq\exp\{r^{\rho_{h}+\varepsilon}\} ~\mbox{ and }~
M(r,\overline{g})=M(r,g)\leq\exp\{r^{\rho_{g}+\varepsilon}\}  ~\mbox{ for any $r>r_{0}$.}
$$
Clearly,
$$M(r,h+\overline{g})\leq M(r,h)+M(r,\overline{g})\leq 2\exp\{r^{\rho_{h}+\varepsilon}\}
$$
so that $\rho\leq\rho_{h}.$ 

Next to prove $\rho\geq \rho_{h}=\max\{\rho_h,\rho_g\}$, we consider the case $\rho_{h}<\infty$.
We need to deal with two subcases.

\bca
$\rho_g<\rho_h$.
\eca
In this case, for any $\varepsilon>0 ~(\varepsilon<\frac{\rho_{h}-\rho_{g}}{3})$,
there exists a sequence of positive real numbers $\{r_{m}\}_{m\geq 1}$ tending to infinity, such that
$$r_{m}>e~\mbox{ and }~ M(r_{m},h)>\exp\{r_{m}^{\rho_{h}-\varepsilon}\} ~\mbox{  for large $m$.}
$$
Also, we take $\{\theta_{m}\}$ such that $|h(r_{m}e^{i\theta_{m}})|=M(r_{m},h)$. Then for large $m$,
\beqq
M(r_{m},h+\overline{g}) &\geq&|h(r_{m}e^{i\theta_{m}})|-|g(r_{m}e^{i\theta_{m}})|\\
&\geq & M(r_{m},h)-M(r_{m},g)\\
&\geq & \exp \big\{r_{m}^{\rho_{h} -\varepsilon}\big\}-\exp \{r_{m}^{\rho_{g}+\varepsilon}\big \}\\
&>& \exp \big \{ r_{m}^{\rho_{h}-2\varepsilon}\big \}
\eeqq
which gives $\rho\geq\rho_{h}-2\varepsilon$. By the arbitrariness of $\varepsilon$, we obtain $\rho\geq \rho_{h}$
and by Definition \ref{de2.2}, i.e., \eqref{eq2.4}, we have $\tau=\tau_{h}$.

\bca
 $\rho_h=\rho_g$.
\eca
In this case, as
$$\frac{n\log n}{\log[1/|a_n|]}\leq \frac{n\log n}{\log[1/(|a_{n}|+|b_{n}|)]},
$$
by Theorem \ref{thm2.1}, we get $\rho_h=\rho_g\leq \rho$ which in turn implies that  $\rho=\rho_h=\rho_g$.
Thus, on one hand, we have
\beqq
r^{-\rho}\log M(r,f)&\leq& r^{-\rho}\log \big(M(r,h)+M(r,g)\big)\\
&\leq & \max\{r^{-\rho}\log 2M(r,h), r^{-\rho}\log 2M(r,g)\}
\eeqq
showing that $\tau\leq \max\{\tau_h, \tau_g\}$ by \eqref{eq2.4}. On the other hand, since
$$n|a_{n}|^{\frac{\rho}{n}}\leq n(|a_{n}|+|b_{n}|)^{\frac{\rho}{n}} ~\mbox{ and }~n|b_{n}|^{\frac{\rho}{n}}\leq n(|a_{n}|+|b_{n}|)^{\frac{\rho}{n}},
$$
by Theorem \ref{thm2.2}, we obtain that $ \max\{\tau_h, \tau_g\}\leq \tau$. Hence $\tau=\max\{\tau_h, \tau_g\}$.

Finally, if $\rho_{g}\leq\rho_{h}=\infty$, it follows similarly that $\rho=\infty$. To sum up, we get $\rho=\rho_{h}$
and the proof is complete.
\epf

\section{Harmonic mappings  with univalent derivatives }\label{csw-sec3}

Denote  by $\mathcal {S}_H$ the class of all univalent sense-preserving
harmonic mappings $f=h+\overline{g}$ in $\ID$, where the following power series expansions
of $h$ and $g$ about the origin are given by
\be\label{eq1.1}
h(z)=z+\sum_{n=2}^{\infty}a_n z^n~ \mbox{ and }~
g(z)=\sum_{n=1}^{\infty}b_n z^n, \quad z\in\ID,
\ee
and where we write for convenience $a_0=0$ and $a_1=1$. Let
\be\label{DPQS-thm5-1} \alpha=\sup\{|a_2|:\, f\in \mathcal {S}_H\}~ \mbox{ and }~\beta=\sup\{|b_2|:\, f\in \mathcal {S}_H\}.\ee
Also, let
$$\mathcal {S}_H^0=\{f\in \mathcal {S}_H:\,f_{\bar{z}}(0)=0\},
$$
so that $\mathcal {S}_H^0\subset\mathcal {S}_H$, and let
$\mathcal{S}=\{f=h+\overline{g}\in \mathcal{S}_H^0:\, g(z)\equiv
0\}.$ {Similar to the central role that $\mathcal {S}$ plays in the
study of univalent function theory}, $\mathcal {S}_H^0$ plays a
vital role in the study of harmonic univalent mappings (see
\cite{AAP,CS,Du,PonRasi2013}).


Let $\mathcal {E}_H$ denote the family of harmonic mappings $f\in \mathcal {S}_H$ which have the property that for each $n\geq 1$,
the harmonic mapping
$f^{(n)}:=h^{(n)}+\overline{g^{(n)}}$ is also univalent in $\mathbb{D}$. Introduce
$$\mathcal{E} =\{f=h+\overline{g}\in \mathcal {E}_H:\,g(z)\equiv 0~\mbox{ in $\ID$}\},
$$
which is a subclass of $\mathcal {S}$ studied, for example in \cite{Bo,ST3, ST0, ST, Sh, Sh-1989}. Let
\be\label{DPQS-eqex7}
\alpha_0=\sup\{|a_2|:\, f\in \mathcal {E}_H\} ~\mbox{ and }~\beta_0=\sup\{|b_2|:\, f\in \mathcal {E}_H\} .
\ee
Since $\mathcal {E}_H\subset\mathcal {S}_H$, it follows that $\alpha_0$ is finite (cf. \cite{AAP, Du}).
Also, $\beta_0$ is finite (cf. \cite{AAP} and \cite[p.~96]{Du}). Also, for functions $h$ and $g$ as in \eqref{eq1.1},
it is easy to see that
$$h^{(n)}(z)=n!a_n +\sum_{k=1}^{\infty}\frac{(n+k)!}{k!} a_{n+k}z^k ~\mbox{ and }~g^{(n)}(z)=n!b_n +\sum_{k=1}^{\infty}\frac{(n+k)!}{k!} b_{n+k}z^k,
$$
which will be used in the proofs of the theorems below.

In \cite{ST}, authors have shown that functions in $\mathcal{E}$ are entire. A natural question is to ask for a harmonic
analog of this statement.

\bprob
{\em
Are functions in $\mathcal {E}_H$ necessarily harmonic entire mappings?
}
\eprob

Our next result answers this question affirmatively.


\bthm\label{thm3.1}
If $f=h+\overline{g}\in \mathcal {E}_H$, then $f$ is a harmonic entire mapping of exponential type, and for every $r$,
we have
$$M(r,f)\leq \frac{e^{2\gamma r}-1}{\gamma}
$$
where $\gamma=\max\{\alpha_0, \beta_0\}$,  $\alpha_0$ and $\beta_0$ are defined by \eqref{DPQS-eqex7}.
\ethm
\bpf Let $f=h+\overline{g}\in \mathcal {E}_H$, where $h$ and $g$ are given by \eqref{eq1.1}.
As $f$ is sense-preserving, we have $J_{f}(0)=|h'(0)|^{2}-|g'(0)|^{2}=1-|b_1|^2>0$
and thus, $|b_1|<1$; and also $\max\{|a_2|, |b_2|\}\leq \gamma$.
By assumption,  $f^{(n)}=h^{(n)}+\overline{g^{(n)}}$ is also univalent in $\mathbb{D}$ for each $n\geq 1$ and thus, without loss of generality,
we assume that $f^{(n)}$ is sense-preserving in $\ID$ so that $J_{f^{(n)}}(0)=((n+1)!)^2[|a_{n+1}|^{2}-|b_{n+1}|^{2}]>0$ which implies that
$a_{n+1}\neq0$ and $|a_{n+1}|>|b_{n+1}|$. Define $F_n$ in $\mathbb{D}$ by
$$F_n(z)=\frac{h^{(n)}(z)-n!a_n}{(n+1)!a_{n+1}}+\overline{\frac{g^{(n)}(z)-n!b_n}{(n+1)!a_{n+1}}}=H_n(z)+\overline{G_n(z)},
$$
where
$$H_n(z)=z+\frac{(n+2)a_{n+2}}{2a_{n+1}}z^2+\cdots\,\,\mbox{and}\,\,G_n(z)=\frac{b_{n+1}}{a_{n+1}}z+\frac{(n+2)b_{n+2}}{2a_{n+1}}z^2+\cdots.
$$
Then $F_n\in \mathcal {E}_H$ and, by \eqref{DPQS-eqex7} and the definition of $\gamma$, it follows that
\be\label{DPQS-eqex3}
|a_{n+2}|\leq \frac{2 \alpha_0 |a_{n+1}|}{n+2}\leq
\frac{2 \gamma |a_{n+1}|}{n+2} ~\mbox{ and }~|b_{n+2}|\leq\frac{2 \beta_0 |a_{n+1}|}{n+2}\leq \frac{2 \gamma |a_{n+1}|}{n+2}.
\ee
An induction argument yields
\be\label{DPQS-eqex2}
|a_{n}|\leq \frac{(2\gamma)^{n-1}}{n!} ~\mbox{ and }~|b_{n}|\leq \frac{(2\gamma)^{n-1}}{n!} ~\mbox{ for $n\geq 1$}
\ee
from which we obtain that $h$ and $g$ are entire, and hence, $f$ is harmonic entire. Moreover, by \eqref{DPQS-eqex2}, we have
for $|z|=r$,
\beqq
|f(z)|\leq |h(z)|+|g(z)| \leq 2\sum_{n=1}^{\infty}\frac{(2\gamma)^{n-1}}{n!}|z|^n
=\frac{e^{2\gamma r}-1}{\gamma}. 
\eeqq
Hence $f$ is a harmonic entire mapping of exponential type.
\epf

We recall that a univalent mapping $f$ harmonic in $\ID$ is called convex if $f(\ID)$ is a convex domain.
A function $f$ ($f'(0)\neq 0$) analytic in $\ID$ is said to be close-to-convex if
the exterior $D=\IC\backslash f(\ID)$ of $f(\ID)$ can be filled with rays emanating on the boundary $\partial D$
and lying completely in $D$. Every close-to-convex function
is known to be univalent in $\ID$. However, close-to-convexity of harmonic mappings requires univalence in
the definition. See \cite{Du} and also \cite{CS,PonRasi2013} for characterizations and discussion on several other
subclasses of ${\mathcal S}_H$.

\beg\label{exam1}
Let $f=h+\lambda\overline{h}$, where $h(z)=\prod_{n=1}^{\infty}(1-\frac{z}{x_{n}})$ is an entire function and
$\lambda$ is a complex constant with $|\lambda|<1$. Following the discussion in \cite[Theorem 1]{Sa},  we let
$$\phi(n)=\log n+1 ~\mbox{ and }~ b(n)=1+\frac{2}{n} ~\mbox{  for $n\geq 1$}.
$$
Also, let $a=e^{2}$ 
so that $a^{\phi(n)}=e^{2}n^{2}$ for $n\geq 1$. Following the second example in
\cite[Section 3]{Sa},   we assume that
$$
|x_{n+1}|\geq e^{2}n^{2}|x_{n}| ~\mbox{ for $n\geq 1$}
$$
and have that for $n\geq 2$,
$$
c(n+1)=1+\frac{1}{e^2 (n+1)^2}+\frac{1}{e^4(n+1)^2(n+2)^2}+ \cdots <\frac{1}{1-[1/(e^2 (n+1)^2)]}\leq 2
$$
and
$$
\sum_{i=2}^{\infty}\frac{1}{e^{2(i-1)}[(i-1)!]^2-2}<\frac{1}{3}.
$$
Let $\ds\big \{x_{j}^{(k)}\big \}_{j=1}^{\infty}$ denote the zeros of $h^{(k)}(z)$. It follows from \cite[Theorem 1]{Sa}  that, for $k=0,1,2...$,
$$\left \{ \begin{array}{ll}
\ds \frac{\big |x_{n+1}^{(k)}\big |}{1+\frac{2}{n}}\leq \big |x_{n}^{(k+1)}\big|< \big |x_{n+1}^{(k)}\big | & \mbox{ for $n\geq 2$},\\[3mm]
\ds 2\big |x_{1}^{(k)}\big |\leq \big |x_{1}^{(k+1)} \big |\leq \big |x_{2}^{(k)}\big |, & \mbox{ and }\\[3mm]
\ds \big |x_{n+1}^{(k)}\big |\geq e^2 n^2\big |x_{n}^{(k)}\big | & \mbox{ for $n\geq 1$.}
\end{array}
\right .
$$
Let $R_{h,n}$ denote the radius of univalence of $h^{(n)}$ at the origin. By using \cite[Theorem 2]{Sa}, we have
$$R_{h,k}=\big |x_{1}^{(k+1)}\big | ~\mbox{ for every $k\geq 0$,}
$$
which implies that if $|x_{1}|\geq \frac{1}{2}$, then $h$ and $h^{(n)}$ are univalent in $\mathbb{D}$.
So $f=h+\lambda\overline{h}$ and $f^{(n)}=h^{(n)}+\lambda\overline{h^{(n)}}$, being composition univalent mappings $h^{(n)}$
with an affine mapping, are univalent in $\mathbb{D}$. \hfill $\Box$
\eeg

\beg\label{exam2}
Let  $f(z)=e^z-1+\lambda\overline{(e^z-1)}$ with $|\lambda|< 1$. Then $f$ and $f^{(n)}$ are all convex and univalent in $\mathbb{D}$.
Thus, $f\in \mathcal {E}_H$.
\hfill $\Box$
\eeg

\beg\label{exam3}
Consider $f=h+\lambda\overline{h}$ with $|\lambda|< 1$, where
$$h(z)=\frac{2}{z}(e^z-(1+z))=\sum_{k=1}^{\infty}a_k z^k, \quad a_k =\frac{2}{(k+1)!}.
$$
We see that
$$H_n(z)=\frac{h^{(n)}(z)-n!a_n}{(n+1)!a_{n+1}}=\sum_{k=1}^{\infty}A_k z^k,
$$
where
$$A_k=\frac{(k+n)!a_{k+n}}{k!(n+1)!a_{n+1}} = \frac{n+2}{k!(k+n+1)}.
$$
As $\{kA_k\}_{k\geq 1}$ is a decreasing sequence, it follows that
$$\sum_{k= 1}^{\infty} |kA_k-(k+1)A_{k+1}|=1 -\lim_{k\rightarrow \infty}(k+1)A_{k+1}=1
$$
and thus, by the result of Ozaki \cite{O1} (see also \cite{Po1}), we conclude that for each $n\geq 0$,
the functions $H_n$ and hence $h^{(n)}$, are close-to-convex (univalent) in $\ID$.  By \cite[5.15. Lemma]{CS}, it follows
that $f$ and $f^{(n)}$ are close-to-convex in $\mathbb{D}$, and hence, $f\in \mathcal {E}_H$.
\hfill $\Box$
\eeg

Let $J_{\nu}$ and $H_{\nu}$ stand for the Bessel and Struve functions of the first kind, and $s_{\mu,\nu}$ be the Lommel function of the first kind
defined in terms of hypergeometric series by \cite{Ba}
$$ s_{\mu,\nu}(z)=\frac{z^{\mu+1}}{(\mu-\nu+1)(\mu+\nu+1)}\cdot {}_{1}F_2\left (1;\frac{\mu-\nu+3}{2},\frac{\mu+\nu+3}{2};-\frac{z^2}{4}\right),
 $$
where $\mu\pm \nu$ are not negative odd integers, and it is a particular solution of the following inhomogeneous
Bessel differential equation
$$z^2w''(z) +zw'(z)+ (z^2-\nu^2)w(z) =z^{\mu+1}.
$$
Here ${}_{1}F_2(a;b,c;z)$ denotes the hypergeometric series given by
$${}_{1}F_2(a;b,c;z):=\sum_{n=0}^{\infty} \frac{(a,n) }{(b,n)(c,n)}\frac{z^n}{n!} =1+\frac{a}{bc}z +\frac{a(a+1)}{b(b+1)c(c+1)}\frac{z^2}{n!} +\cdots,
$$
where $(a,0)=1$ for $a\neq 0$ and, for each positive integer $n$, $(a,n):=a(a+1)\cdots (a+n-1)$.

\beg\label{exam4}
Consider $f=h+\lambda\overline{h}$ with $|\lambda|< 1$,  and $h$ is one of $f_{\nu}$, $h_{\nu}$ and $l_{\nu}$, where
$$
f_{\nu}(z)=2^{\nu}\Gamma(\nu+1)z^{1-\frac{\nu}{2}}J_{\nu}(\sqrt{z})=\sum_{n=0}^{\infty}\frac{(-1)^n\Gamma(\nu+1)z^{n+1}}{4^n n!\Gamma(\nu+n+1)}
$$
with $\nu\geq \nu_0$ ($\nu_0\simeq-0.5623...$ is the unique root of equation $f_{\nu}'(1)=0$ on $(-1,\infty)$),
$$
h_v(z)=\sqrt{\pi}2^{\nu}\Gamma(\nu+\frac{3}{2})z^{\frac{1-\nu}{2}}H_{\nu}(\sqrt{z})
=\frac{\sqrt{\pi}}{2}\sum_{n=0}^{\infty}
\frac{(-1)^n\Gamma(\nu+\frac{3}{2})z^{n+1}}{4^n\Gamma(n+\frac {3}{2})\Gamma(\nu+n+\frac {3}{2})}, ~|\nu|\leq \frac{1}{2},
$$
and
\beqq
l_{\mu}(z)&=& \ds \mu(\mu+1)z^{-\frac{\mu}{2}+\frac{3}{4}}s_{\mu-\frac{1}{2},\frac{1}{2}}(\sqrt{z})
=z \cdot {}_{1}F_2\left (1;\frac{\mu+2}{2},\frac{\mu+3}{2};-\frac{z}{4}\right )\\
&=&\sum_{n=0}^{\infty}\frac{(-1)^n \Gamma(\frac{\mu}{2}+1)\Gamma(\frac{\mu}{2}+\frac{3}{2})
z^{n+1}}{4^n \Gamma(\frac{\mu}{2}+n+1)\Gamma(\frac{\mu}{2}+n+\frac{3}{2})}
\eeqq
with $\mu\in (-1,1)$ and $\mu\neq 0$.
By using \cite[Theorems 1, 2, 3]{Ba}, we
obtain that $h$ is starlike, and all $h^{(n)}$ ($n\geq 1$) are close-to-convex in $\ID$.
This implies that $f$ and $f^{(n)}$ are close-to-convex, and then $f\in \mathcal {E}_H$.
\hfill $\Box$
\eeg

\br
{\em
There are ways to establish entire functions $h$ such that $h$ and all its derivatives $h^{(n)}$
are convex (univalent) in the unit disk $\ID$. Such functions can be used to generate harmonic entire mappings
$f$ and $f^{(n)}$ such that they are close-to-convex in $\mathbb{D}$.
}
\er

From Theorem \ref{thm3.1}, we know that if  $f^{(n)}=h^{(n)}+\overline{g^{(n)}}$
$(n\geq 0)$ are all univalent in $\mathbb{D}$, then $f$ must be a harmonic entire mapping of exponential type.
In fact, to obtain that $f$ is harmonic entire, it is not necessary that all $f^{(n)}=h^{(n)}+\overline{g^{(n)}}$
$(n\geq 1)$ are univalent in $\mathbb{D}$. For a given fixed $n$, let $R_n$ be the largest positive number with the
property that $f^{(n)}$ is univalent in $|z|<R_n$. Note that $R_n$ is finite unless $f^{(n)}$ is a polynomial
(i.e. both $h^{(n)}$ and $g^{(n)}$ are some polynomials). Next, we consider the relationship between the sequence
$\{R_n\}_{n=1}^{\infty}$ and the radius of convergence {$R_f$} of $f$ about the origin. First we recall
the following lemma, see \cite[Theorem 1]{AZ} or \cite[Theorem 3.1]{QW}.

\blem\label{thm3.2-lem} Let $f(z)=z+\sum_{k=2}^{\infty}a_k z^k+\overline{\sum_{k=1}^{\infty}b_kz^k}$ be harmonic in $\mathbb{D}$, $|b_1|<1$ and
$$
\sum_{n=2}^{\infty}n(|a_n|+|b_n|)\leq 1-|b_1|.
$$
Then $f$ is sense-preserving and univalent in $\mathbb{D}$.
\elem

In fact in \cite[Lemma 1]{KPV-14},   the conclusion of Lemma \ref{thm3.2-lem} strengthened by concluding that $ |f_z(z)-1|<1-|f_{\overline{z}}(z)|$ in $\mathbb{D}$
and $f$ is close-to-convex in $\mathbb{D}$.

\bthm\label{thm3.2}
Let $f(z)=\sum_{k=1}^{\infty}a_k z^k+\overline{\sum_{k=1}^{\infty}b_kz^k}$ belong to $\mathcal{S}_{H}$
 Also, suppose that
$R_f$ denotes the radius of convergence of $f$ and $R_n$ denotes the radius of univalence of $f^{(n)}$. Then,
\be\label{thm3.2-1}
\liminf_{n\rightarrow \infty} n R_n\leq 2 (\alpha+\beta) R_f,
\ee
where $\alpha$ and $\beta$ are given by \eqref{DPQS-thm5-1}.

If $\limsup_{n\rightarrow \infty }\big|\frac{b_n}{a_n}\big|=\delta<1$ or
$\limsup_{n\rightarrow \infty }\big|\frac{a_n}{b_n}\big|=\delta<1$, then
$$\frac{1-\delta}{1+\delta}R_f\log 2\leq \limsup_{n\rightarrow \infty} n R_n.
$$
If there is a nonnegative integer $N$ such that for $n\geq N$, $R_n>0$, then
$$\liminf_{n\rightarrow\infty}[n(R_N R_{N+1}\cdots R_{n})^{1/n}]\leq 2\max\{\alpha, \beta\}e R_f.
$$

If $\limsup_{n\rightarrow \infty }\big|\frac{b_n}{a_n}\big|=\delta<1$ or
$\limsup_{n\rightarrow \infty }\big|\frac{a_n}{b_n}\big|=\delta<1$, and  there is a positive integer $N$, such that $\{A_{n-1}\}_{n=N}^{\infty}$, $A_{n-1}=\frac{|a_{n-1}|+|b_{n-1}|}{|a_n|+|b_n|}$,
is a positive and nondecreasing sequence, then
$$\frac{1-\delta}{1+\delta}R_f\log 2\leq \liminf_{n\rightarrow\infty}n R_n\leq \limsup_{n\rightarrow\infty}n R_n\leq 2 (\alpha+\beta) R_f.
$$
\ethm\bpf
If an infinite number of the $R_n$ are zero, then \eqref{thm3.2-1} is obviously true.
So suppose that there is a nonnegative integer $N$ such that $R_n>0$ for all $n\geq N$. Note that this
implies that $|a_{n+1}|\neq|b_{n+1}|$ for $n\geq N$, because $J_{f^{(n)}}(0)\neq 0$. Now, we
consider $f^{(n)}(z)=h^{(n)}(z)+\overline{g^{(n)}(z)}$, where
$$
h^{(n)}(z) =n!a_n +\sum_{k=1}^{\infty}\frac{(n+k)!}{k!} a_{n+k}z^k
$$
and a similar expression holds for $g^{(n)}(z)$. Without loss of generality,
we assume that $f^{(n)}$ is sense-preserving in $|z|< R_n$. In particular, $J_{f^{(n)}}(0)>0$
which gives the condition $|a_{n+1}|>|b_{n+1}|$ and thus,
$$\frac{R_n^{-1}h^{(n)}(R_n z) -n!a_n}{(n+1)!a_{n+1}} +\overline{\frac{R_n^{-1}g^{(n)}(R_n z) -n!b_n}{(n+1)!a_{n+1}}}
\in \mathcal{S}_{H}.
$$
Then as in \eqref{DPQS-eqex3} we see that
\be\label{thm3.2-5}
(n+2)|a_{n+2}|R_n\leq 2\alpha |a_{n+1}| ~\mbox{ and }~  (n+2)|b_{n+2}|R_n \leq 2\beta |a_{n+1}|.
\ee
Summing up the two inequalities in \eqref{thm3.2-5} gives
\be\label{thm3.2-3}
(n+2)R_n\leq 2(\alpha+\beta)\frac{|a_{n+1}|}{|a_{n+2}|+|b_{n+2}|}  \leq 2(\alpha+\beta)A_{n+1}.
\ee
Obviously, the series $\sum_{k=1}^{\infty}(|a_k|+|b_k|)z^k$ converges for $|z|<R_f$ which shows that
$$\liminf_{n\rightarrow \infty}\frac{|a_{n+1}|+ |b_{n+1}|}{|a_{n+2}|+|b_{n+2}|}\, =\liminf_{n\rightarrow \infty} A_{n+1}\leq R_f,$$
and then \eqref{thm3.2-3} shows that
$$\liminf_{n\rightarrow \infty} n R_n\leq 2 (\alpha+\beta) R_f.
$$
Moreover, from \eqref{thm3.2-5}, we find that
\be\label{thm3.2-4}
|a_{n+2}|\leq |a_{n+1}| \frac{2\alpha}{(n+2)R_n}  ~\mbox{ and }~ |b_{n+2}| \leq |a_{n+1}| \frac{2\beta}{(n+2)R_n}.
\ee
Using \eqref{thm3.2-4}, an induction argument shows that
$$|a_n|\leq \frac{2^{n-N-1}\alpha^{n-N-1} |a_{N+1}| (N+1)!}{(R_N R_{N+1}\cdots R_{n-2})n!} ~\mbox{ for $n\geq N+2$}
$$
and a similar expression holds (with $\beta$ in place of $\alpha$) for the second term in \eqref{thm3.2-4}. In particular,
\be\label{thm3.2-7}
\max\{|a_n|, |b_n|\}\leq \frac{2^{n-N-1}{\gamma_1}^{n-N-1} |a_{N+1}| (N+1)!}{(R_N R_{N+1}\cdots R_{n-2})n!} ~\mbox{ for $n\geq N+2$},
\ee
where $\gamma_1=\max\{\alpha, \beta\}$.

Now we recall that there are two positive numbers, $A=\sqrt{2\pi}$ and $B=\sqrt{2\pi}e^{1/24}$ such that (cf. \cite[p.~183]{Mi})
\be\label{thm3.2-6}
A n^{1/2}(n/e)^n < n! < B n^{1/2}(n/e)^n ~\mbox{ for $n\in \IN$}.
\ee
Using this and \eqref{thm3.2-7}, we have
$$
\frac{1}{R_f}=\limsup_{n\rightarrow\infty}\max\{|a_n|^{1/n}, |b_n|^{1/n}\}\leq
\frac{2\gamma_1 e}{\liminf_{n\rightarrow\infty}[n(R_N R_{N+1}\cdots R_{n})^{1/n}]},
$$
that is,
$$
\liminf_{n\rightarrow\infty}[n(R_N R_{N+1}\cdots R_{n})^{1/n}]\leq 2\gamma_1 e R_f.
$$
By assumption, there is a positive integer $N$ such that $\{A_{n-1}\}_{n=N}^{\infty}$
is a positive and nondecreasing sequence. Note that this implies that $R_n>0$. Hence
$$\lim_{n\rightarrow \infty}A_{n}= R_f,$$
which by using \eqref{thm3.2-3} yields
$$\limsup_{n\rightarrow \infty} n R_n\leq 2(\alpha+\beta)R_f.
$$

Let $0<r<R_f$. Since $\sum |a_n|r^n<\infty$ and $\sum |b_n|r^n<\infty$ which imply that there is an increasing sequence
$\{n_p\}_{p=1}^{\infty}$ of positive integers such that for $p=1, 2, \ldots$, and $k=2, 3, \ldots$,
$$
|a_{n_p+1}|+|b_{n_p+1}|\geq (|a_{n_p+k}|+|b_{n_p+k}|)r^{k-1}.
$$
Let $x_n=n(1-2^{-1/(n+2)})$. Then

\vspace{7pt}

$\ds\sum_{k=2}^{\infty}\frac{(n_p+k)!(|a_{n_p+k}|+|b_{n_p+k}|)r^k (x_{n_p})^k}{(k-1)!{n_p}^k}
$
\beqq
&\leq&(|a_{n_p+1}|+|b_{n_p+1}|)r \sum_{k=2}^{\infty}\frac{(k+n_p)!}{(k-1)!}\left(\frac{x_{n_p}}{n_p}\right)^k\\
&\leq&\frac{(n_p+1)!(|a_{n_p+1}|+|b_{n_p+1}|)r x_{n_p}}{n_p} \sum_{k=2}^{\infty}\frac{\big((n_p+2)(1-2^{-1/(n_p+2)})\big)^{k-1}}{(k-1)!}\\
&=&\frac{(n_p+1)!(|a_{n_p+1}|+|b_{n_p+1}|)r x_{n_p}}{n_p} \left(e^{(n_p+2)(1-2^{-1/(n_p+2)})}-1\right)\\
&\leq& \frac{(n_p+1)!(|a_{n_p+1}|+|b_{n_p+1}|)r x_{n_p}}{n_p},
\eeqq
so that

\vspace{7pt}

$\ds \sum_{k=2}^{\infty}\frac{(k+n_p)!(|a_{n_p+k}|+|b_{n_p+k}|)r^k (x_{n_p})^k}{(k-1)!n_p^k} B_{n_p}^{k}
$
\beqq
&\leq& \sum_{k=2}^{\infty}\frac{(k+n_p)!(|a_{n_p+k}|+|b_{n_p+k}|)r^k (x_{n_p})^k}{(k-1)!n_p^k} B_{n_p}^{2}\\
&\leq&\frac{(n_p+1)!\big| |a_{n_p+1}|-|b_{n_p+1}|\big|r x_{n_p}B_{n_p}}{n_p},
\eeqq
where
$$
B_{n_p}= \frac{\big| |a_{n_p+1}|-|b_{n_p+1}|\big|}{|a_{n_p+1}|+|b_{n_p+1}|}.
$$

Define $F_{p}$ in $\mathbb{D}$ by $$
F_{p}(z)=f^{(n_p)}\Big (\frac{rzx_{n_p}B_{n_p}}{n_p}\Big ).
$$
By using Lemma \ref{thm3.2-lem}, we find that $F_p$ is univalent in $\mathbb{D}$.
Hence if $\limsup_{n\rightarrow \infty }\big|\frac{b_n}{a_n}\big|=\delta<1$ or
$\limsup_{n\rightarrow \infty }\big|\frac{a_n}{b_n}\big|=\delta<1$, we have
$$r x_{n_p}B_{n_p}\leq n_p R_{n_p},$$
and then
$$\frac{1-\delta}{1+\delta}R_f\log 2 \leq \limsup_{n\rightarrow \infty} n R_n.
$$
If $f$ is such that for some integer $N$, $\{A_{n-1}\}_{n=N}^{\infty}$
is a positive and nondecreasing sequence, and $\limsup_{n\rightarrow \infty }\big|\frac{b_n}{a_n}\big|=\delta<1$ or
$\limsup_{n\rightarrow \infty }\big|\frac{a_n}{b_n}\big|=\delta<1$.
If $n\geq N$, we let $r_n=A_{n-1}$ so that
$$|a_{n+1}|+|b_{n+1}|\geq (|a_{n+k}|+|b_{n+k}|)r_{n}^{k-1} ~\mbox{ for $k\geq 2$}.
$$
Using the preceding arguments, it follows that
$$ r_{n} x_{n}B_n\leq n R_{n} \mbox{ for $n\geq N$},
$$
which, since $\lim_{n\rightarrow\infty}r_nx_n=R_f\log 2$,  shows that
$$
\frac{1-\delta}{1+\delta}R_f\log 2\leq \liminf_{n\rightarrow\infty}n R_n.
$$
The proof is complete.
\epf

 For the inequality \eqref{thm3.2-1}, we suspect that there is a better estimate of $R_n$. Here we give it as a conjecture.

\begin{conj}\label{thm3.2-conj}
Let $f(z)=\sum_{k=1}^{\infty}a_k z^k+\overline{\sum_{k=1}^{\infty}b_kz^k}\in \mathcal{S}_{H}$.
 Also, suppose that
$R_f$ denotes the radius of convergence of $f$ and $R_n$ denotes the radius of univalence of $f^{(n)}$. Then
$$
\liminf_{n\rightarrow \infty} n R_n\leq 2 \max\{\alpha,\beta\} R_f.
$$
\end{conj}

We say that a harmonic entire mapping $f=h+\overline{g}$ is called a harmonic polynomial if both
$h$ and $g$ are some polynomials. Thus, we say that $f$ is a transcendental harmonic entire mapping
if $f$ is not a harmonic polynomial.

\bcor\label{thm3.2-cor1}
Let $f(z)=\sum_{k=0}^{\infty}a_k z^k+\overline{\sum_{k=1}^{\infty}b_kz^k}$,  and $R_n$ denote the radius of univalence of $f^{(n)}=h^{(n)}+\overline{g^{(n)}}$.
We have the following:
\bee
\item If $\lim_{n\rightarrow\infty}nR_n=\infty$, then $f$ is a transcendental
harmonic entire mapping.
\item If $f$ is a transcendental harmonic entire mapping, and either
$\limsup_{n\rightarrow \infty }\big|\frac{b_n}{a_n}\big|=\delta<1$ or
$\limsup_{n\rightarrow \infty }\big|\frac{a_n}{b_n}\big|=\delta<1$ holds, then $\limsup_{n\rightarrow\infty}nR_n=\infty$.
\item Suppose that either
$\limsup_{n\rightarrow \infty }\big|\frac{b_n}{a_n}\big|=\delta<1$ or
$\limsup_{n\rightarrow \infty }\big|\frac{a_n}{b_n}\big|=\delta<1$ holds. If $f$ is a transcendental harmonic entire mapping, and
there is a positive integer $N$, such that $\Big \{\frac{|a_{n-1}|+|b_{n-1}|}{|a_n|+|b_n|}\Big \}_{n=N}^{\infty}$
is a positive and nondecreasing sequence, then $\lim_{n\rightarrow\infty}n R_n=\infty$.
\eee
\ecor

Theorem \ref{thm3.2} is a generalization of \cite[Theorem 1]{ST0}. From  \eqref{thm3.2-1}, we know that if
$R_n$ converges to zero slowly enough, then $f$ must still be harmonic entire. From \cite[p.~315]{ST},  it is known that the
converse of Corollary \ref{thm3.2-cor1} is false.

  In the following, we consider relationships between the growth of
$\{R_n\}_{n=0}^{\infty}$ and  the order of $f$.

\bthm\label{thm3.3} Suppose that $f(z)=h(z)+\overline{g(z)}=\sum_{k=0}^{\infty}a_k z^k+\overline{\sum_{k=1}^{\infty}b_kz^k}$
is a transcendental harmonic entire mapping of order $\rho$,
and $R_n$ is the radius of univalence of $f^{(n)}=h^{(n)}+\overline{g^{(n)}}$. Then
\be\label{thm3.3-1}
\liminf_{n\rightarrow\infty}\frac{\log(\max\{1, n R_n\})}{\log n}\leq \frac{1}{\rho}.
\ee
\ethm
\bpf To prove \eqref{thm3.3-1}, we may assume that $n R_n\geq 1$ for all $n$ and that $\rho>0$.
Let $\rho_h$ and $\rho_g$ be the orders of $h$ and $g$, respectively. We may
assume that $\rho_h\geq\rho_g$.

If $\rho_h=\rho_g$ then, by Theorem \ref{thm2.3}, we have $\rho=\rho_h=\rho_g$. Without lose of generality, we suppose that there is an increasing sequence
$\{n_p\}_{p=1}^{\infty}$ of positive integers such that, for $p=1,2,\cdots$, $f^{(n_p)}$ is sense-preserving in $|z|<R_{n_p}$(otherwise we consider $\overline{f}$).
It is known from \cite{Sh-1946} that
$$\frac{1}{\rho}=\frac{1}{\rho_h}\geq \liminf_{n\rightarrow\infty}\frac{\log|a_n/a_{n+1}|}{\log n}.
$$
It follows from \eqref{thm3.2-5} that
$$
\log n_p R_{n_p}\leq \log (2\alpha)+\log |a_{n_p+1}/a_{n_p+2}|,
$$
and hence,
$$\frac{1}{\rho}\geq \liminf_{p\rightarrow\infty}\frac{\log n_p}{\log(n_p+1)}\left(\frac{\log n_p R_{n_p}}{\log n_p}-\frac{ \log 2\alpha}{\log n_p}\right)=\liminf_{n\rightarrow\infty}\frac{\log nR_n}{\log n}.
$$

Now the consider the case $\rho_h>\rho_g$. Then $\rho=\rho_h$ by Theorem \ref{thm2.3}.
 We claim that there is an increasing sequence $\{n_p\}_{p=1}^{\infty}$ of positive integers such that, for $p=1,2,\cdots$, $f^{(n_p)}$ is sense-preserving in $|z|<R_{n_p}$.
Otherwise, there is an integer $N_0>0$ such that for $n>N_0$, $f^{(n)}$ is sense-reversing in $|z|<R_{n}$, and then
$|a_n|<|b_n|$, which implies that the maximum term of $h$ is less than the maximum term of $g$. By \cite[2.3.5]{Bo1}, $\rho_h\leq\rho_g$ which contradicts
the assumption that $\rho_h>\rho_g$. By using similar arguments as that in the case $\rho_h=\rho_g$ we obtain that
$$\frac{1}{\rho}\geq \liminf_{n\rightarrow\infty}\frac{\log nR_n}{\log n}.
$$

\epf

Theorem \ref{thm3.3} is a generalization of \cite[Theorem 2]{ST} to the case of harmonic mappings.

It is known that if $f(z)=z+\sum_{n=2}^{\infty}a_nz^n\in\mathcal{S}$, then
$|a_n|\leq n$ for $n\ge 2$ (see \cite{dB85}). If $f=h+\overline{g}\in \mathcal {S}_H^0$, where
$h$ and $g$ are given by \eqref{eq1.1} with $b_1=0$, then
the corresponding coefficient conjecture due to Clunie and Sheil-Small \cite{CS} is
$$|a_n|\leq \frac{(2n+1)(n+1)}{6}~\mbox{ and }~|b_n|\leq \frac{(2n-1)(n-1)}{6} ~ \mbox{for $n\ge 2$}.
$$
However, only the elementary inequality $|b_2|\leq 1/2$ has been verified so for, thanks to the Schwarz lemma.
For the best known bound for $|a_2|$, we refer to \cite{AAP}. On the other hand, the above conjecture has been
verified for harmonic starlike mappings, harmonic close-to-convex mappings and for typically real harmonic mappings,
respectively (cf. \cite{CS,Du}).  The analog coefficient inequalities for the classes of harmonic convex mappings are
established by Clunie and Sheil-Small \cite{CS}.

Moreover, for $f=h+\overline{g}\in \mathcal {S}_H$,  where $h$ and $g$ are given by \eqref{eq1.1}, it is conjectured
that \cite{CS}
\be\label{eq3.1}
|a_n|<\frac{2n^2+1}{3} ~\mbox{ and }~ |b_n|<\frac{2n^2+1}{3}  ~ \mbox{for $n\ge 2$}.
\ee

In the following, we will consider the family $\mathcal{S}_{H}^{c}$ of harmonic entire mappings $f=h+\overline{g}$ such that
$f\in \mathcal {S}_H$ with the power series expansions given by \eqref{eq1.1} and satisfying the
condition \eqref{eq3.1}.

%

\bthm\label{thm3.4}
Let $f(z)=h(z)+\overline{g(z)}=\sum_{k=0}^{\infty}a_k z^k+\overline{\sum_{k=1}^{\infty}b_kz^k}$
be harmonic  in $\mathbb{D}$. Let $\{n_p\}_{p=1}^{\infty}$ be a strictly increasing sequence of nonnegative
integers such that
$$\frac{f^{(n_p)}-n_p!a_{n_p}-n_p!\overline{b_{n_p}}}{(n_p+1)! a_{n_p+1}}\in \mathcal{S}_{H}^{c}\; \mbox{if} \;|a_{n_p+1}|\geq |b_{n_p+1}|,$$
and
$$\frac{\overline{f^{(n_p)}}-n_p!\overline{a_{n_p}}-n_p!b_{n_p}}{(n_p+1)! b_{n_p+1}}\in \mathcal{S}_{H}^{c}\; \mbox{if} \;|a_{n_p+1}|\leq |b_{n_p+1}|.$$
Let $\lambda=\liminf_{p\rightarrow \infty}\frac{n_p}{n_{p+1}}$. Then we have the following:
\bee
\item If $\lambda=1$, then $R_f=\infty$, i.e. $f$ is entire.
\item If $0<\lambda<1$, then $R_f\geq \lambda^{\frac{\lambda}{\lambda-1}}/(1-\lambda)$.
\item If $\lambda=0$, then $R_f \geq 1$.
\eee
(While Part {\rm (3)} is obvious, it is included here for completeness.)
\ethm
\bpf
From the assumption, we known that $f^{(n_p)}=h^{(n_p)}+\overline{g^{(n_p)}}$ is univalent in $\mathbb{D}$. 
 Without loss of generality, we may assume that there is a subsequence of $\{n_p\}_{p=1}^{\infty}$ (we still use the notation $\{n_p\}_{p=1}^{\infty}$) such that, for $p=1,2, \cdots$, $f^{(n_p)}$ are sense-preserving
(otherwise, we may consider $\overline{f^{(n_p)}}$). Then  $J_{f^{(n_p)}}(0)>0$ which implies that
$a_{n_p+1}\neq0$ and $|a_{n_p+1}|>|b_{n_p+1}|$. As
$$H_{n_p}(z)=\frac{h^{(n_p)}(z)-n_p!a_{n_p}}{(n_p+1)! a_{n_p+1}}=  z + \frac{1}{(n_p+1)! a_{n_p+1}}\sum_{k=2}^{\infty}\frac{(n_p+k)!}{k!} a_{n_p+k}z^k
$$
and
$$
G_{n_p}(z)=\frac{g^{(n_p)}(z)-n_p! b_{n_p}}{(n_p+1)! \overline{a_{n_p+1}} } =\frac{b_{n_p+1}}{\overline{a_{n_p+1}}}z+ \frac{1}{(n_p+1)! \overline{a_{n_p+1}}}\sum_{k=2}^{\infty}\frac{(n_p+k)!}{k!} b_{n_p+k}z^k,
$$
we see that $F_{n_p} =H_{n_p}+\overline{G_{n_p}} \in \mathcal{S}_{H}^{c}$, by assumption. Using \eqref{eq3.1}, we have, for
$p=1, 2, \ldots$ and $k=2, 3, \ldots$,
\be\label{DPQS-eqex5}
|a_{n_p+k}|\leq \frac{2k^2+1}{3}\frac{k!(n_p+1)!}{(n_p+k)!}|a_{n_p+1}|.
\ee
%
Setting $k=n_{p+1}-n_p+1$ and applying the method of induction on $p$, we obtain, for $p=2, 3, \ldots$,
$$|a_{n_p+1}| \leq \frac{(\frac{1}{3})^{p-1}(n_1+1)! |a_{n_1+1}| }{(n_p+1)!}
\times \prod_{j=1}^{p-1}(2(n_{j+1}-n_j+1)^2+1)(n_{j+1}-n_j+1)!.
$$
By \eqref{DPQS-eqex5}, we conclude that, for $p\geq 2$ and $2\leq k\leq n_{p+1}-n_p+1$,
\beq\label{thm3.4-1}
 \nonumber |a_{n_p+k}|  &\leq& \frac{(\frac{1}{3})^{p}(n_1+1)! |a_{n_1+1}| (2k^2+1)k!}{(n_p+k)!}\\
\nonumber &&\times \prod_{j=1}^{p-1}(2(n_{j+1}-n_j+1)^2+1)(n_{j+1}-n_j+1)!\\
  &<&\frac{(n_1+1)! |a_{n_1+1}|  k^2 k!}{(n_p+k)!}
\times \prod_{j=1}^{p-1}(n_{j+1}-n_j+1)^2(n_{j+1}-n_j+1)!.
\eeq
A similar expression holds for $|b_{n_p+k}|$.

If $f^{(n_p)}$ is sense-reversing then \eqref{thm3.4-1} holds with $|b_{n_1+1}|$
on the right in place of $|a_{n_1+1}|$. More generally, for $p\geq 2$ and $2\leq k\leq n_{p+1}-n_p+1$, we have
\beq\label{thm3.4-1a}
\max\{|a_{n_p+k}|, |b_{n_p+k}|\}  &<&\frac{(n_1+1)! \max\{|a_{n_1+1}|, |b_{n_1+1}|\}  k^2 k!}{(n_p+k)!}
\\&&
 \nonumber \times \prod_{j=1}^{p-1}(n_{j+1}-n_j+1)^2(n_{j+1}-n_j+1)!.
\eeq


Using \eqref{thm3.2-6}  on the right-hand side of \eqref{thm3.4-1}, taking the $(n_p+k)$-th
root on both sides of the resulting inequality, and applying \cite[Lemma 2]{ST} to part of
right-hand-side of this, it follows that for $p\geq 2$ and for $2\leq k\leq n_{p+1}-n_p+1$,

\beq
\nonumber |a_{n_p+k}|^{\frac{1}{n_p+k}}&<&\left [ (n_1+1)!\max \{|a_{n_1+1}|, |b_{n_1+1}|\}B^p e^{n_1+1-p} \sqrt{\frac{k^5}{n_p+k}}\, \right]^{1/(n_p+k)}\\
\nonumber &&\times \left(1+\frac{n_p}{p}\right)^{7p/(2n_p)}\frac{k^{k/(n_p+k)}}{n_p+k}
\prod_{j=1}^{p-1}(n_{j+1}-n_j+1)^{\frac{n_{j+1}-n_j}{n_p+k}}.
\eeq
By using the similar arguments as that in the proof of \cite[Theorem]{ST}, it follows that
\beq\label{thm3.4-5}
\nonumber |a_{n_p+k}|^{\frac{1}{n_p+k}}&<&\left [ (n_1+1)!\max \{|a_{n_1+1}|, |b_{n_1+1}|\} B^p e^{n_1+1-p} \sqrt{\frac{k^5}{n_p+k}}\; \right]^{1/(n_p+k)}\\
\nonumber&& \times \left(1+\frac{n_p}{p}\right)^{7p/(2n_p)}  \max\left \{\frac{4^{1/n_p}}{n_p}\prod_{j=1}^{p-1}(n_{j+1}-n_j+1)^{\frac{n_{j+1}-n_j}{n_p}},\right.\\
&&\left.\hspace{0.3cm}\,\,\frac{n_{p+1}^{1/(n_{p+1})}}{n_{p+1}}\prod_{j=1}^{p}(n_{j+1}-n_j+1)^{\frac{n_{j+1}-n_j}{n_{p+1}}}  \right\}.
\eeq
Letting $\gamma_0=\limsup_{p\rightarrow\infty} p/n_p $,  we have
\beq\label{thm3.4-3}
\nonumber \frac{1}{R_f}&\leq& \limsup_{k\rightarrow\infty}|a_{k}|^{1/k}\\
\nonumber &=&\limsup\left \{|a_{n_p+k}|^{\frac{1}{n_p+k}}: \,2\leq p, ~ 2\leq k\leq n_{p+1}-n_p+1\right \}\\
&\leq & K \limsup_{p\rightarrow\infty}\frac{1}{n_p}\prod_{j=1}^{p-1}(n_{j+1}-n_j+1)^{\frac{n_{j+1}-n_j}{n_p}},
\eeq
where
$$K= \left\{
\begin{array}{cl}
1                                                          &\mbox{ if $\gamma_0=0$}\\
\ds \left[\frac{B}{e}\left (1+\frac{1}{\gamma_0}\right )^{7/2}\right]^{\gamma_0}&\mbox{ if $\gamma_0>0$}
\end{array} \right.  .
$$

By using the similar arguments as that in the proof of \cite[Theorem]{ST} again, it follows from \eqref{thm3.4-3} that
$$
\frac{1}{R_f}\leq \left\{
\begin{array}{cl}
0  &\mbox{ if $\lambda=1$}\\
\ds K(1-\lambda)\lambda^{\lambda/(1-\lambda)}    &\mbox{ if $0<\lambda<1$}\\
\ds K                                  &\mbox { if $\lambda=0$}
\end{array} \right. .
$$

From the assumption, $R_f\geq 1$. Then for the case $\lambda=0$, it is obvious $R_f\geq 1$.
Suppose $0<\lambda<1$. Let $0<\varepsilon<\lambda$. Then \cite[Lemma]{ST} shows that there is
 a subsequence of $\{n_p\}_{p=1}^{\infty}$, $\{n_{p_k}\}_{k=1}^{\infty}$, such that
 $\lambda-\varepsilon\leq \liminf_{k\rightarrow\infty}\frac{n_{p_k}}{n_{p_{k+1}}}$,
$\limsup_{k\rightarrow\infty}\frac{n_{p_k}}{n_{p_{k+1}}}<1$ and $\lim_{k\rightarrow\infty}k/n_{p_k}=0$. Applying what we have just proved to this subsequence
$\{n_{p_k}\}_{k=1}^{\infty}$, we have $K=1$, and then
$$
\frac{1}{R}\leq (1-\lambda+\varepsilon)(\lambda-\varepsilon)^{(\lambda-\varepsilon)/(1-\lambda+\varepsilon)}.
$$
Since this is true for all sufficiently small $\varepsilon$, we have
$$
\frac{1}{R}\leq (1-\lambda)\lambda^{\lambda/(1-\lambda)}.
$$
The proof is complete.
\epf

Our theorem is stated so as to get information about $R_f$ when $\lambda$ is known.
 If $R_f$ is known, then Theorem \ref{thm3.4} gives a bound on $\lambda$.

\bcor\label{thm3.4-cor1}
Let $f$ and $\{n_p\}_{p=1}^{\infty}$ be defined as in Theorem {\rm \ref{thm3.4}}. If $f$ cannot be extended to a harmonic mapping in a
disk $|z|<R_f$ with the radius $R_f>1$, then
$$\liminf_{p\rightarrow \infty}\frac{n_{p-1}}{n_p}=0.
$$
\ecor

If $\lim_{p\rightarrow \infty}\frac{n_{p-1}}{n_p}=1$, then Theorem \ref{thm3.4} shows that $f$ is harmonic entire.
From this, it follows that, if $n_p-n_{p-1}=o(n_p)$, then also $f$ is harmonic entire. In fact, to insure that $f$
is harmonic entire, it is enough to put this condition on the second differences of the  $n_p$.

\bcor\label{thm3.4-cor2}
Let $f$ and $\{n_p\}_{p=1}^{\infty}$ be defined as in Theorem \ref{thm3.4}. If
$$n_{p+2}-2n_{p+1}+n_p=o(n_p),
$$
then $f$ is a harmonic entire mapping.
\ecor

In the following, we investigate relationship between the order of a harmonic entire mapping $f$ and the number of its derivatives $f^{(n)}$
which are univalent in a given disk.

\bthm\label{thm3.5}
Let $f(z)=h(z)+\overline{g(z)}=\sum_{k=0}^{\infty}a_k z^k+\overline{\sum_{k=1}^{\infty}b_kz^k}$ be a harmonic entire mapping
of order $\rho$. Let $\{n_p\}_{p=1}^{\infty}$ be a strictly increasing sequence of nonnegative integers
such that $$\frac{f^{(n_p)}-n_p!a_{n_p}-n_p!\overline{b_{n_p}}}{(n_p+1)! a_{n_p+1}}\in \mathcal{S}_{H}^{c}\; \mbox{if} \;|a_{n_p+1}|\geq |b_{n_p+1}|,$$
$$\frac{\overline{f^{(n_p)}}-n_p!\overline{a_{n_p}}-n_p!b_{n_p}}{(n_p+1)! b_{n_p+1}}\in \mathcal{S}_{H}^{c}\; \mbox{if} \;|a_{n_p+1}|\leq |b_{n_p+1}|.$$
Then
\be\label{DPQS-eqex6}
\rho\leq \frac{1}{\ds 1-\limsup_{p\rightarrow\infty}\frac{\log (n_p-n_{p-1})}{\log n_p}}.
\ee
\ethm
\bpf By using the similar arguments as that at the beginning of the proof of Theorem \ref{thm3.4}, without loss of generality, we assume that
$f^{(n_p)}$ are sense-preserving for $p=1,2, \cdots$. Then
$$\frac{h^{(n_p)}-n_p!a_{n_p}}{(n_p+1)!a_{n_p+1}}+\overline{\frac{g^{(n_p)}-n_p! b_{n_p}}{(n_p+1)!\overline{a_{n_p+1}}}} \in \mathcal{S}_{H}^{c}.
$$
From \cite[2.4.1 Theorem]{Bo1}, the order of $h$ is same as the order of  $h^{(n_p)}$ and the order of $g$ is same
as the order of $g^{(n_p)}$. Moreover, Theorem \ref{thm2.3} shows that the order of $f$ is same as
the order of $f^{(n_p)}$. Since $f^{(n_p)}$ is sense-preserving, we have that $h^{(n_p)}$ is an entire function of order $\rho$.
It follows after certain calculations (see the proof of Theorem \ref{thm3.4}, especially inequality \eqref{thm3.4-5}) that
there is a positive number $K^{*}$ such that, for $p\geq 2$ and $2\leq k\leq n_{p+1}-n_p+1$,
\beq
\nonumber \max\left \{|a_{n_p+k}|^{\frac{1}{n_p+k}}, |b_{n_p+k}|^{\frac{1}{n_p+k}}\right \}&\leq& K^{*} \max\left\{\frac{1}{n_p}\prod_{j=1}^{p-1}(n_{j+1}-n_j+1)^{\frac{n_{j+1}-n_j}{n_p}},\right.\\
\nonumber &&\,\,\,\,\,\,\, \left. \frac{1}{n_{p+1}}\prod_{j=1}^{p}(n_{j+1}-n_j+1)^{\frac{n_{j+1}-n_j}{n_{p+1}}}  \right\}.
\eeq
 By using the similar arguments as that of the proof of  \cite[Theorem 1]{Sh}, we obtain the desired
inequality \eqref{DPQS-eqex6}.
\epf

We observe that Theorem \ref{thm3.5} exhibits an upper bound on the order $\rho$ of $f$ in terms of the
sequence $\{n_p\}_{p=1}^{\infty}$.

\bcor\label{thm3.5-cor1}
Let $f$, $\rho$ and $\{n_p\}_{p=1}^{\infty}$ be as in Theorem {\rm \ref{thm3.5}}. If
$$
\log(n_p-n_{p-1})=o(\log n_p),
$$
then $\rho \leq 1$. If $n_p-n_{p-1}\leq \mu$ for all large $p$, then $f$ is of
exponential type no bigger than $\sqrt{2\pi}e^{-47/24}(\mu+1)^{9/2}$.
\ecor
\bpf The first part of the conclusion follows directly from Theorem \ref{thm3.5}. To prove the second part, we may assume that
$n_p-n_{p-1}\leq \mu$ for all $p$. Then by \eqref{thm3.4-1a} we have, for $p\geq 2$ and $2\leq k\leq n_{p+1}-n_p+1$, that
$$\max \{|a_{n_p+k}|,
|b_{n_p+k}|\} < \frac{D^{*} k^2 k!}{(n_p+k)!}
\prod_{j=1}^{p-1}(n_{j+1}-n_j+1)^2(n_{j+1}-n_j+1)!,
$$
where $D^{*}$ is a constant. Now, by using \eqref{thm3.2-6}, the last inequality gives
\beqq
&&\hspace{-3cm} \ds \max \left \{(n_p+k)|a_{n_p+k}|^{1/(n_p+k)},
(n_p+k)|b_{n_p+k}|^{1/(n_p+k)}\right \}\\
 &\leq& \left ( \frac{D^{*}B^pe^{n_1+1-p}}{A}\left (\frac{k^5}{n_p+k}\right )^{1/2}\right)^{1/(n_p+k)} \\
&&\times k^{k/(n_p+k)}
 \prod_{j=1}^{p-1}(n_{j+1}-n_j+1)^{(n_{j+1}-n_j+7/2)/(n_{p}+k)}.
\eeqq
By using the similar arguments as that of the proof of \cite[Corollary 1]{ST}, we have
\beqq
\limsup_{k\rightarrow\infty}k|a_k|^{1/k} &=&  \limsup \{(n_p+k)|a_{n_p+k}|^{1/(n_p+k)}:\, p\geq 2,\, 2\leq k\leq n_p-n_{p-1}+1\}\\
 &\leq&\frac{B}{e}(\mu+1)^{9/2}=e\Big(\sqrt{2\pi}e^{-47/24}(\mu+1)^{9/2}\Big).
\eeqq
Similarly,
$$\limsup_{k\rightarrow\infty}k|b_k|^{1/k}
\leq\frac{B}{e}(\mu+1)^{9/2}=e\Big(\sqrt{2\pi}e^{-47/24}(\mu+1)^{9/2}\Big).
$$
By \cite[2.2.10 Theorem]{Bo1},   we see that the type of $h$ and $g$ are all less than $\sqrt{2\pi}e^{-47/24}(\mu+1)^{9/2}$.
The second part of the corollary follows from Theorem \ref{thm2.3}.
\epf

\bthm\label{thm3.6}
Let $f=h+\overline{g}$ be harmonic in $\mathbb{D}$ and $\{n_p\}_{p=1}^{\infty}$ be as in Theorem {\rm \ref{thm3.5}.} If
$\rho^{*}$, defined by
$$\frac{1}{\rho^{*}}=1-\limsup_{p\rightarrow\infty}\frac{\log
(n_p-n_{p-1})}{\log n_p},
$$
is finite, then $f$ is a harmonic entire mapping of order no greater than $\rho^{*}$. If $\rho^{*}=\infty$, then $f$ need not be harmonic
entire, and if it is harmonic entire, then it may be of any order.
\ethm
\bpf
Suppose that $\rho^{*}<\infty$. Then $\lim_{p\rightarrow\infty}(n_p/n_{p+1})=1$ and Theorem \ref{thm3.6}
implies that the hypotheses of Theorem \ref{thm3.4} are satisfied showing that $f$ is harmonic entire of order no
greater than $\rho^{*}$.

The functions exhibited in Theorem 6 of \cite{ST3} show that it may be the case that $\rho^{*}=\infty$
and yet the radius of convergence is $1$. This establishes the second part.
\epf

\br
{\em
Theorem {\rm \ref{thm3.6}} shows that Theorem {\rm \ref{thm3.5}} is sharp.
}
\er

\subsection*{Acknowledgements}
The third author is supported by NSF of Hebei Science Foundation (No. A2018201033).
The work of the second author was supported by Mathematical Research Impact Centric Support (MATRICS) of
the Department of Science and Technology (DST), India  (MTR/2017/000367).

%

\end{document}